# Diophantine equation of degree sixteen

Oliver Couto & Seiji Tomita

## Abstract

While there is not much publications, about degree sixteen Diophantine equation we do have an identity given by Ramanujan (ref. #1). Also on the internet even though there are numerical solutions to degree sixteen for eg. (16-7-24) equation (ref. #5) there are hardly any parametric solutions. An Octic degree parameterization has been arrived at by Choudhry & Zagar (ref. 2). The authors have given a parametric solution to the equation: $(a^4 - -b^4)(c^4 - d^4)(e^8 - f^8) = (u^4 - -v^4)(w^4 - x^4)(y^8 - z^8)$. We have also given numerical solution but because of the high degree (sixteen) we only get a minimum integer value for the variables at more than five digits. We have also given some new identities related to degree four & eight.

Consider the below equation:

$$(a^4 - -b^4)(c^4 - d^4)(e^8 - f^8) = (u^4 - v^4)(w^4 - x^4)(y^8 - z^8) ---(1)$$

Equation (1) was derived from using the below quartic equation:

$$(a^4 - b^4)(c^4 - d^4) = m(p^4 - q^4)$$

## Section (A)

Consider the below equation:

$$m(a^4 - b^4) = (c^4 - d^4)(e^4 - f^4)$$

Theorem,

$$m(a^4 - b^4) = (c^4 - d^4)(e^4 - f^4) \quad \text{---- (A)}$$

Eqn (A) has infinitely many integer solutions,
where, $m = n(p^2 + q^2)$ and n,p,q are arbitrary integer.

Proof,

We split eqn (A) into two simultaneous eqns.

$$(p^2 + q^2)(a^2 + b^2) - (c^2 + d^2)(e^2 + f^2) = 0 \quad \ldots \ldots \ldots \ldots (1)$$

$$n(a^2 - b^2) - (c^2 - d^2)(e^2 - f^2) = 0 \quad \ldots \ldots \ldots \ldots \ldots (2)$$

From eq(1), let $\{a, b, c, d\} = \{fu - ev,\ eu + fv,\ pu + qv, qu - pv\}$

Then eq(2) becomes to

$$(nv^2 - nu^2 - p^2 * u^2 - 4puqv - q^2v^2 + q^2u^2 + p^2v^2)e^2$$
$$-4nfuev + nf^2u^2 - p^2f^2v^2 + p^2f^2u^2 - nf^2v^2 + q^2f^2v^2 + 4puqvf^2 - q^2f^2u^2 \ldots \ldots (3)$$

For the quadratic in eq(3) to have rational solutions, the discriminant must be a rational square, then we get,

$$w^2 = (-2nq^2 + 2np^2 - 2p^2q^2 + n^2 + q^4 + p^4)u^4$$
$$+ (8npq - 8pq^3 + 8p^3q)vu^3$$
$$+ (2n^2 - 2q^4 - 2p^4 - 4np^2 + 4nq^2 + 20p^2q^2)v^2u^2$$
$$+ (-8p^3q + 8pq^3 - 8npq)v^3u$$
$$+ (-2nq^2 + 2np^2 - 2p^2q^2 + n^2 + q^4 + p^4)v^4.$$

Let, $U = \frac{u}{v}, W = \frac{w}{v^2},$ then we get quartic eqn:

$W^2 =$

$$(p^2 - q^2 + n)^2 U^4 + (8pqn - 8pq^3 + 8p^3q)U^3$$
$$+ (2n^2 + 4q^2n + 20p^2q^2 - 4p^2n - 2q^4 - 2p^4)U^2 + (8pq^3 - 8pqn - 8p^3q)U$$
$$+ (p^2 - q^2 + n)^2 \ldots \ldots \ldots \ldots \ldots (4)$$

This quartic has a rational point,

$$Q(U, W) = (0, p^2 - q^2 + n),$$

so is birationally equivalent to an elliptic curve below.

$$Y^2 - 8pqYX + (32p^3qn - 32p^3q^3 + 16p^5q - 32q^3 * pn + 16q^5p + 16pqn^2)Y$$
$$= X^3 + (2n^2 - 4p^2n + 4q^2n - 2q^4 - 2p^4 + 4p^2q^2)X^2$$
$$+ (16p^2q^6 + 48p^4q^2n - 48p^2q^4n + 48p^2q^2n^2 + 16q^2n^3 - 24q^4n^2 - 16p^2n^3 - 16p^6n$$
$$- 24p^4q^4 + 16q^6n + 16p^6q^2 - 24p^4n^2 - 4p^8 - 4q^8 - 4n^4)X$$
$$+ 104q^8n^2 - 96q^6n^3 - 160p^6q^6 - 48p^2q^{10} + 120p^4q^8 - 48p^{10}q^2 - 48q^{10}n + 120p^8q^4$$
$$+ 48p^{10}n + 104p^8n^2 + 96p^6n^3 + 24p^4n^4 - 16p^2n^5 + 24q^4n^4$$
$$+ 16q^2n^5 + 8p^{12} + 8q^{12} - 8n^6 - 416p^2q^6n^2 + 240p^2q^8n - 480p^4q^6n - 288p^4q^2n^3$$
$$+ 624p^4q^4n^2 + 480p^6q^4n - 416p^6q^2n^2 - 240p^8q^2n + 288p^2q^4n^3$$
$$- 48p^2q^2n^4$$

This elliptic curve has a point

$$P(X, Y) = (-2n^2 + 4p^2n - 4q^2n + 2q^4 + 2p^4 - 4p^2q^2, -32pqn^2).$$

According to Nagell-Lutz theorem, this point P is of infinite order, and the multiples kP, k = 2, 3, ...give infinitely many points.
Then simultaneous equations (1),(2) has infinitely many integer solutions.

2Q(U) corresponding to 2P(X,Y) is

$$U = \frac{4qp(p^2 - q^2 + n)}{(q-p)(q+p)(-q^2 + p^2 + 2n)}$$

then we get,

$$a = p^8 + 4p^7q + (-4q^2 + 4n)p^6 + (8nq - 12q^3)p^5 + (6q^4 + 4n^2 + 4nq^2)p^4 +$$
$$(-16nq^3 + 12q^5)p^3 + (8q^2n^2 - 4q^4n - 4q^6)p^2 + (8nq^5 - 4q^7)p +$$
$$(4n^2q^4 - 4nq^6 + q^8)$$

$$b = p^8 - 4p^7q + (-4q^2 + 4n)p^6 + (-8nq + 12q^3)p^5 + (6q^4 + 4n^2 + 4nq^2)p^4$$
$$+ (16nq^3 - 12q^5)p^3 + (8q^2n^2 - 4q^4n - 4q^6)p^2 + (-8nq^5 + 4q^7)p$$
$$+ 4n^2q^4 - 4nq^6 + q^8$$

$$c = 3p^4q + (-2q^3 + 2nq)p^2 + 2nq^3 - q^5$$

$$d = p^5 + (2q^2 + 2n)p^3 + (-3q^4 + 2nq^2)p$$

$$e = p^4 + (-2q^2 + 2n)p^2 - 4npq - 2nq^2 + q^4$$

$$f = p^4 + (-2q^2 + 2n)p^2 + 4npq - 2nq^2 + q^4$$

We substitute, $n = q^2$ & we get,

$$(a, b, c, d, e, f) = (p^4 + 4pq^3 - q^4)^2, (p^4 - 4pq^3 - q^4)^2,$$
$$(p^5 + 4p^3q^2 - pq^4, 3p^4q + q^5), (p^4 + 2p^3q - 2p^2q^2 + 2pq^3 + q^4), (p^4 - 2p^3q - 2p^2q^2 - 2pq^3 + q^4)$$

---

Also, (a,b), becomes an eight power.

So we now have an eqn shown below:

$$m(u^8 - v^8) = (c^4 - d^4)(e^4 - f^4) \quad \text{-------(a)}$$

$$\text{where, } (a, b) = (u^2, v^2)$$

$$m = (n)(p^2 + q^2) = q^2(p^2 + q^2) \text{ since we have } n = q^2$$

## Section (B)

Consider the below eqn:

$$(a^4 - b^4)(c^4 - d^4) = m_1(y^8 - z^8) \quad \text{--- (1)}$$

$$(u^4 - v^4)(w^4 - x^4) = m_2(e^8 - f^8) \; ----(2)$$

where,

$$m_1 = q^2(p^2 + q^2)$$

$$m_2 = s^2(r^2 + s^2)$$

And eqn (1) & (2) are of the form parametrized in section (A) above as eqn (a)

we take, $\quad m_1 = m_2$

Dividing eqn (1) by eqn (2) & cross multiplying we get:

$$(a^4 - b^4)(c^4 - d^4)(e^8 - f^8) = (u^4 - v^4)(w^4 - x^4)(y^8 - z^8)$$

Now since we need, $\quad m_1 = m_2$

$$m_1 = q^2(p^2 + q^2)$$

$$m_2 = s^2(r^2 + s^2)$$

Hence, $\quad q^2(p^2 + q^2) = s^2(r^2 + s^2)$

Above is parametrized as:

$$p = (2k^2 + 22k - 7)$$

$$q = 8(k + 1)(k - 2)$$

$$r = 2(8k^2 - 2k + 17)$$

$$s = 4(k^2 - k - 2)$$

For, $k = 1$ we get:

$$(p, q, r, s) = (17, 16, 46, 8)$$

From section (A) we have parametric form for,

(a,b,c,d,e,f) & (u,v,w,x,y,z) as below:

$$a = p^4 + 4pq^3 - q^4$$

$$b = p^4 - 4pq^3 - q^4$$

$$c = p^5 + 4p^3q^2 - pq^4$$

$$d = 3p^4q + q^5$$

$$e = r^4 + 2r^3s - 2r^2s^2 + 2rs^3 + s^4$$

$$f = r^4 - 2r^3s - 2r^2s^2 - 2rs^3 + s^4$$

---

$$u = r^4 + 4rs^3 - s^4$$

$$v = r^4 - 4rs^3 - s^4$$

$$w = r^5 + 4r^3s^2 - rs^4$$

$$x = 3r^4s + s^5$$

$$y = p^4 + 2p^3q - 2p^2q^2 + 2pq^3 + q^4$$

$$z = p^4 - 2p^3q - 2p^2q^2 - 2pq^3 + q^4$$

For, (p,q,r,s)=( 17,16,46,8), we get for:

$$(a^4 - b^4)(c^4 - d^4)(e^8 - f^8) = (u^4 - v^4)(w^4 - x^4)(y^8 - z^8)$$

where:

a = 296513, b = 260543, c = 5336657, d = 5057584, e = 5815184, f = 2606224

u = 4567568, v = 4379152, w = 230692576, x = 107491712, y = 297569, z = 295391

A different eight degree parametric Identity is given below:
$$(p^4 - q^4) = m(r^4 - s^4)(t^4 - u^4) --- (1)$$
where:

$$p = 2(50m^2 - 37m + 5)$$

$$q = 2(10m^2 + 13m - 5)$$

$$r = 3(3m - 1)$$

$$s = 7m - 1$$

$$t = 10m - 1$$

$$u = 10m - 7$$

for m = 3 we get:

For m = 3 we get:

$$(86^4 - 31^4) = 3(6^4 - 5^4)(29^4 - 23^4)$$

---

A sixteen degree parametric identity is given below:
Consider the below eqn:

$$(m - n)(u - v) = 4(x - y)(z - w) ------ (1)$$

where:
$$m = a^8 + b^8 + c^8$$
$$n = d^8 + e^8 + f^8$$
$$u = p^8 + q^8 + r^8$$
$$v = s^8 + t^8 + u^8$$

$$x = (ab)^4 + (bc)^4 + (ca)^4$$
$$y = (de)^4 + (df)^4 + (ef)^4$$
$$z = (pq)^4 + (pr)^4 + (qr)^4$$
$$w = (st)^4 + (su)^4 + (tu)^4$$

$$(m - n)(u - v) = 4(x - y)(z - w)$$

Hence we have:

$$[(a^8 + b^8 + c^8) - (d^8 + e^8 + f^8)]$$
$$= 2[((ab)^4 + (bc)^4 + (ca)^4) - ((de)^4 + (df)^4 + (ef)^4)] \quad --- (2)$$
$$[(p^8 + q^8 + r^8) - (s^8 + t^8 + u^8)]$$
$$= 2[((pq)^4 + (pr)^4 + (qr)^4) - ((st)^4 + (su)^4 + (tu)^4)] \quad --- (3)$$

Above eqn (2) &(3) is satisfied at:
$$(a, b, c, d, e, f) = (732, 804, 342, 293, 513, 536)$$

$$(p, q, r, s, t, u)) = (63232, 71825, 76032, 104593, 61776, 88400)$$

And $(a, b, c, d, e, f)$ & $(p, q, r, s, t, u)$ are having parametrization as given in ref. (2) by Choudhry & Zargar paper

On multiplying eqn (2) &(3) above, we get eqn (1) above which is a sixteen degree equation.

---

We have the below quartic equation:

$$a^4 + b^4 + ab(a^2 + ab + b^2) = c^4 + d^4 + cd(c^2 + cd + d^2) \quad --- (1)$$

Taking:

$$[a = pt + m, b = qt + n, c = pt - m, d = qt - n] \quad --- (2)$$

In the above eqn, (1), we noticed there was a pattern in the numerical solutions for
(p,q) to eqn. (1)

we took: $\quad p = u^2 \ \& \ q = v^2$

After substituting above (and using maple software) we then took,
$$(m, n) = ((u^2 - 2uv + 2v^2), (2u^2 - 2uv + v^2))$$
and we noticed that eqn (1) is satisfied at $t = [(u + v/(u - v)]$.

Thus we get the below parametric solution:

$$a = u^3 - u^2v + 2uv^2 - v^3$$

$$b = v^3 + 2u^2v - uv^2 - u^3$$

$$c = u^3 - 2u^2v + 2u*v^2$$

$$d = 2vu^2 - 2uv^2 + v^3$$

For, $(u, v) = (8,7)$ we get numerical solution:

$$(a, b, c, d) = (101, 67, 91, 80)$$

Ramanujan equation:

Ramanujan, gave a sixteen degree parametric identity & is shown below:

$$45((a^8 + b^8 + (a+b)^8) - (d^8 + e^8 + (d+e)^8))^2 =$$

$$64((a^{10} + b^{10} + (a+b)^{10}) - (d^8 + e^8 + (d+e)^8)) *$$

$$((a^6 + b^6 + (a+b)^6) - (d^6 + e^6 + (d+e)^6))$$

Condition is: $a^2 + ab + b^2 = c^2 + cd + d^2$ ----- (1)

Above eqn (1) has parameterization:

$$(a, b, c, d) = ((x + y + 1), (xy - 1), (xy + y + 1), (x - y))$$

$$for, (x, y) = (3,2) \text{ we get: } (a, b, c, d) = (6,5,9,1)$$

Hence the degree sixteen numerical solution is,

$$45((6^8 + 5^8 + 11^8) - (9^8 + 1^8 + 10^8))^2 =$$

$$64((6^{10} + 5^{10} + 11^{10}) - (9^{10} + 1^{10} + 10^{10})) * ((6^6 + 5^6 + 11^6) - (9^6 + 1^6 + 10^6))$$

-------------------------------------------

### References:

1) Degree sixteen equation, Ramanujan, Wolfram math world, Diophantine equation

   https://mathworld.wolfram.com/Ramanujan 6-10-8.

2. Octic Diopantine equation. A. Choudhry & A. Zargar,

   Arxiv, 2206.14084v1 dated  28 June 2022.

3. Generalized Parametric Solutions to Multiple Sums of Powers

   Oliver Couto & Seiji Tomita, Universal Journal of Applied

   Mathematics  3(5): 102-111, 2015  http://www.hrpub.org

   DOI:  10.13189/ujam.2015.030503

4. A new diophantine equation involving fifth powers, Oliver Couto &

   Ajai Choudhry, Acta Arithmetica journal,  15 June 2021,

   DOI: 10.4064/aa210419-4-7.

5. Website, Oliver Couto, www.celebrating-mathematics.com

6. Jaroslaw Wroblewski & Ajai Choudhry. Ideal solutions of the Tarry-Escott

   problem of degree eleven with applications to sums of thirteenth powers.

   Hardy-Ramanujan Journal, Hardy-Ramanujan Society, 2008, 31, pp.1-13.

7. Website, Seiji Tomita website: www.maroon.dti.ne.jp/fermat/eindex.html

8. Tito Piezas web site for degree eleven and

   higher powers,(sites.google.com/site/tpiezas/031)

9. Database numerical solutions degree 16,

   etc.  (http://euler.free.fr/database.txt)

Contact info:

Oliver Couto

University of Waterloo,

Ontario, Canada.

Email: matt345@celebrating-mathematics.com

Seiji Tomita

Japan Computer Company

Tokyo, Japan.

Email:fermat@m15.alpha-net.ne.jp

---